\documentclass[12pt]{amsart}
\usepackage{amssymb,amsfonts,amsmath}
\usepackage[top=1in,left=1in,bottom=1in,right=1in]{geometry}
\usepackage{graphicx}
\usepackage[usenames,dvipsnames]{color}
\usepackage{enumerate}
\newtheorem{theorem}{Theorem}[section]
\newtheorem{proposition}[theorem]{Proposition}
\newtheorem{example}[theorem]{Example}
\newtheorem{lemma}[theorem]{Lemma}

\renewcommand{\proof}{\noindent {\it Proof: }}
\renewcommand{\qed}{\null\hfill$\Box$\medskip}
\newcommand{\R}{\mathbb{R}}

\newcommand{\Z}{\mathbb{Z}}
\newcommand{\Q}{\mathbb{Q}}
\newcommand{\C}{\mathbb{C}}

\begin{document}
\title{On Volumes of Permutation Polytopes}
\author{Katherine Burggraf, 
Jes\'{u}s De Loera, 
Mohamed Omar}
\date{\today}
\maketitle


\begin{abstract}
This paper focuses on determining the volumes of permutation polytopes associated to cyclic groups, dihedral groups, groups of automorphisms of tree graphs, and Frobenius groups. We do this through the use of triangulations and the calculation of Ehrhart polynomials. We also present results on the theta body hierarchy of various permutation polytopes.
\end{abstract}

\section{Introduction}~\label{Introduction}

Volumes are fundamental geometric invariants of convex bodies. However,
 volumes of lattice polytopes often have additional algebraic and
combinatorial meaning. For example, \emph{order polytopes} have
volumes equal to the number of linear extensions of their associated
poset \cite{Stanley:TwoPosetPolytopes:1986}. The Catalan number $C_n$ appears as the volume of the
convex hull of the positive root configuration ${A}_{n}^{(+)}$
\cite{Gelfand+Graev+Postnikov:PositiveRoots:1997}. Catalan numbers also appear as factors of the volume of the Chan-Robbins-Yuen polytope
\cite{ChanRobbinsYuen2000,Zeilberger1999}. Volumes of lattice polytopes are also interesting in the
context of algebraic geometry, as they can be used to determine the number of
solutions of zero-dimensional systems of polynomial equations and the degrees of algebraic varieties
\cite{Bernsteinetal1976,Sturmfels1996,Sturmfels2002}. This motivates the investigation of volumes of
polytopes arising as convex hulls of finite groups of permutation
matrices.

Perhaps the key example of such a polytope is the
\emph{Birkhoff polytope} $B_n$.  It is defined as the convex hull of
all $n \times n$ permutation matrices, or equivalently, as the convex
hull of the natural permutation representation of the symmetric group
$S_n$; see \cite{BeckPixton2003, Brualdi1988, BrualdiGibson1977,
  CanfieldMcKay2009, ChanRobbins1999, DeLoeraLiuYoshida2009,
  DiaconisGangolli1995, Pak2000, Sturmfels1996} and references therein
for a summary of its known properties. Subpolytopes of $B_n$ have
been shown to have remarkably beautiful properties; see
\cite{Ahmed2008, Brualdi1988, BrualdiLiu1991, ChanRobbins1999,
  ChanRobbinsYuen2000, DiaconisWood2010, Mirsky1961,
  SchreckTinhofer1988, Zeilberger1999} and references therein. 
This is particularly true for permutation polytopes, those polytopes that
 arise by taking convex hulls of permutation representations of special subgroups of $S_n$
with concrete sets of generators. Their geometry reflects their
group-theoretic structure. For known results on permutation
polytopes, see \cite{BaumeisterHaaseNillPaffenholz2009,
  BrualdiLiu1991, CollinsPerkinson2008, CunninghamWang2004,
  GuralnickPerkinson2006, HoodPerkinson2004, Steinkamp1999} and
references therein. For definitions pertaining to permutation polytopes, see \cite{BaumeisterHaaseNillPaffenholz2009,GuralnickPerkinson2006}.

Before stating our results, we will clarify some terminology. The \emph{normalized
volume} of a $d$-dimensional polytope $P \subset \R^n$ with respect to an affine lattice $L \subset \R^n$ is the volume form that assigns a volume
of one to the smallest $d$-dimensional simplices in $\R^n$ whose
vertices are in $L$. The \emph{volume} of $P$ is
its normalized volume in the lattice aff$(P) \cap \Z^n$.  We say $P$
is \emph{unimodular} with respect to $L$ if it has a triangulation whose simplices are all unimodular; that is, the vertices of
any simplex in the triangulation span the lattice $L$.  For more
on triangulations with respect to particular lattices used in this paper,
see Section~\ref{Preliminaries}.
In what follows, we identify the symmetric
group $S_n$ on $\{1,2,\ldots,n\}$ through its representation by $n
\times n$ permutation matrices; that is, for any $g \in S_n$, we
identify $g$ with the $n \times n$ matrix whose $(i,j)$-entry is one if $g(i)=j$ and 0 otherwise. We denote the identity by $e$
throughout.  We denote a subgroup $G$ of $S_n$ by $G \leq S_n$.  Such
a subgroup is called a \emph{permutation group}.  For any permutation
group $G \leq S_n$, we refer to the polytope $P(G) := \mbox{conv}\{g
\ | \ g \in G\}$ as the \emph{permutation
  polytope} associated to $G$.

\begin{example}\label{example1}
Let $G \leq S_4$ be the group consisting of the four permutations $\{e,(1 \ 2),(3 \ 4),(1 \ 2)(3 \ 4)\}$.  Then $P(G)$ is the convex hull of the matrices
\[
\begin{pmatrix}
1 & 0 & 0 & 0 \\
0 & 1 & 0 & 0 \\
0 & 0 & 1 & 0 \\
0 & 0 & 0 & 1
\end{pmatrix},
\ \
\begin{pmatrix}
0 & 1 & 0 & 0 \\
1 & 0 & 0 & 0 \\
0 & 0 & 1 & 0 \\
0 & 0 & 0 & 1
\end{pmatrix},
\ \
\begin{pmatrix}
1 & 0 & 0 & 0 \\
0 & 1 & 0 & 0 \\
0 & 0 & 0 & 1 \\
0 & 0 & 1 & 0
\end{pmatrix},
\ \
\begin{pmatrix}
0 & 1 & 0 & 0 \\
1 & 0 & 0 & 0 \\
0 & 0 & 0 & 1 \\
0 & 0 & 1 & 0
\end{pmatrix}.
\]
This polytope is geometrically a square.  Now let $H \leq S_4$ be the group consisting of the four permutations $\{e,(1 \ 2)(3 \ 4),(1 \ 3)(2 \ 4),(1 \ 4)(2 \ 3)\}$.  Then $P(H)$ is the convex hull of the matrices
\[
\begin{pmatrix}
1 & 0 & 0 & 0 \\
0 & 1 & 0 & 0 \\
0 & 0 & 1 & 0 \\
0 & 0 & 0 & 1
\end{pmatrix},
\ \
\begin{pmatrix}
0 & 1 & 0 & 0 \\
1 & 0 & 0 & 0 \\
0 & 0 & 0 & 1 \\
0 & 0 & 1 & 0
\end{pmatrix},
\ \
\begin{pmatrix}
0 & 0 & 1 & 0 \\
0 & 0 & 0 & 1 \\
1 & 0 & 0 & 0 \\
0 & 1 & 0 & 0
\end{pmatrix},
\ \
\begin{pmatrix}
0 & 0 & 0 & 1 \\
0 & 0 & 1 & 0 \\
0 & 1 & 0 & 0 \\
1 & 0 & 0 & 0
\end{pmatrix}.
\]
This polytope is geometrically a tetrahedron.
\end{example}
\normalsize

\vspace{0.1in}

Note that Example~\ref{example1} shows that the geometric structure of a permutation polytope depends
 on the presentation of the group that defines it (i.e., on the
choice of generators).  Both of the examples above are groups isomorphic to $(\Z/2\Z)^2$, but
their permutation polytopes are not even combinatorially isomorphic.

Our focus in this paper is on determining the volumes, or normalized volumes, of permutation
polytopes associated to cyclic groups, dihedral groups, groups of
automorphisms of tree graphs, and Frobenius groups. We have two essential approaches to computing the volume of a polytope. One approach
is to triangulate the polytope by simplices of equal volume. 
Another approach is to obtain the volume through
the Ehrhart polynomial of the polytope, a polynomial which counts integer points in dilations of the polytope and whose leading coefficient is the volume of the polytope. For more on these approaches, see Section~\ref{Preliminaries}.  

To begin our study, we introduce our first two classes of permutation polytopes.
The \emph{cyclic group} $C_n \leq S_n$ 
is the group generated by the permutation $(1
\ 2 \ \cdots \ n)$.  The \emph{dihedral group} $D_n \leq S_n$ is the group generated
by the permutations $r = (1 \ 2 \ \cdots \ n)$ and $f = (1 \ n)(2
\ n-1) \cdots (\lfloor \frac{n+1}{2} \rfloor \ \lceil \frac{n+1}{2}
\rceil)$.  In
Section~\ref{cyclicdihedral}, we determine particular unimodular
 triangulations of $P(C_n)$ and $P(D_n)$ with respect to the lattices
 aff$(P(C_n)) \cap \Z^{n \times n}$ and aff$(P(D_n)) \cap \Z^{n \times n}$
 respectively. This allows us to recover their volumes via their Ehrhart polynomials.

\begin{theorem}\label{thm:cyclicdihedralehrhart}
Let $n$ be an integer, $n > 2$.
\begin{enumerate}
\item 
The volume of $P(C_n)$ is $\frac{1}{(n-1)!}$. The Ehrhart polynomial of $P(C_n)$ is $\binom{t+n-1}{n-1}$.

\vspace{0.1in}

\item  If $n$ is odd, the volume of $P(D_n)$ is $\frac{n}{(2n-2)!}$. The Ehrhart polynomial of $P(D_n)$ is
\[
\sum_{k=0}^{n-2} \binom{2n}{k+1} \binom{t-1}{k} + \sum_{k=n-1}^{2n-2} \left( \binom{2n}{k+1} - \binom{n}{k-n+1} \right) \binom{t-1}{k}.
\]

\vspace{0.1in}

\item If $n$ is even, $n=2m$, the volume of $P(D_n)$ is $\frac{n^2}{4 \cdot (2n-3)!}$. The Ehrhart polynomial of $P(D_n)$ is
\begin{align*}
\sum_{k=0}^{m-2} \binom{2n}{k+1} \binom{t-1}{k} - \sum_{k=m-1}^{2m-2} \left( \binom{2n}{k+1} - 2\binom{2n-m}{k+1-m} \right) \binom{t-1}{k} +& \\
\sum_{k=2m-1}^{4m-3} \left( \binom{2n}{k+1} - 2 \binom{2n-m}{k+1-m} + \binom{2n-2m}{k+1-2m} \right) \binom{t-1}{k}.&
\end{align*}
\end{enumerate}
\end{theorem}

In Section~\ref{Frobenius}, we study \emph{Frobenius polytopes} as
introduced by Collins and Perkinson in \cite{CollinsPerkinson2008}.
These are permutation polytopes $P(G)$ where $G$ is a \emph{Frobenius
  group}.  A group
$G \leq S_n$ is Frobenius if it has a proper subgroup $H$ such that for all $x \in G \backslash H$, $H \cap (xHx^{-1}) = \{e\}$.  The special subgroup $H$ is known as the \emph{Frobenius complement} of $G$ and is unique up to conjugation. Moreover, every Frobenius group $G \leq S_n$ has a special proper subgroup $N$ of size $n$ called the \emph{Frobenius kernel} which consists of the identity and all elements of $G$ that have no fixed points; see Chapter 16 of \cite{AlperinBell1995}.  The Frobenius kernel and Frobenius complement have trivial intersection, and $G = NH$.
 The class of Frobenius groups includes semi-direct products of cyclic groups, some matrix groups over finite fields, the alternating group $A_4$, and many others.  See \cite{Wielandt} for more on Frobenius groups.  We determine triangulations of Frobenius polytopes and a formula for their normalized volumes, in particular showing that the normalized volumes are completely characterized by the size of the Frobenius complement and the size of the Frobenius kernel.

\vspace{0.1in}

\begin{theorem}~\label{thm:Frobenius}
Let $G \leq S_n$ be a Frobenius group with Frobenius complement $H$ and Frobenius kernel $N$.  The normalized volume of $P(G)$ in the sublattice of $\Z^{n \times n}$ spanned by its vertices is \[ \frac{1}{(|H||N|-|H|)!} \sum_{\ell = 0}^{\lfloor \frac{|H|(|N|-1)-1}{|N|} \rfloor} \binom{(|H|-\ell)|N|}{(|H|-\ell)|N| - |H|+1} \binom{|H|-1}{\ell} {(-1)}^{\ell}.\]
\end{theorem}

In order to better understand these polytopes, we can approximate them with a sequence of special convex bodies. In the 1980s, L. Lov\'asz approximated the stable set polyhedra from graph theory using a convex body called the \emph{theta body}; see \cite{Lovasz1994}. In \cite{GouveiaParriloThomas2010}, the authors generalize Lov\'asz's theta body for $0/1$ polytopes (that is, polytopes whose vertices have coordinates of zero and one) to generate a sequence of semidefinite programming relaxations of the convex hull of the common zeroes of a set of real polynomials; see \cite{Lovasz1994, LovaszSchrijver1991} and Section~\ref{Preliminaries} for more on this topic.  We briefly study the theta body hierarchy of our permutation polytopes.  For instance, we prove that convergence of the first iterate always occurs for Frobenius groups.  This implies many structural results, such as the existence of reverse lexicographic unimodular triangulations.  See \cite{Sullivant2006} for more on this.

\begin{proposition}\label{prp:frobeniustwolevel}
If $G \leq S_n$ is a Frobenius group, then $P(G)$ is two-level and hence $G$ is $TH_1$-exact.
\end{proposition}

\vspace{0.1in}

In Section~\ref{AutBinTrees}, we develop a method for computing the Ehrhart polynomial of $P(G)$ when $G$ is the automorphism group of a rooted binary tree on $n$ vertices.  This method relates the Ehrhart polynomials of permutation polytopes associated to direct products and wreath products of groups to the Ehrhart polynomials of the individual permutation polytopes themselves.  A key theorem in this regard is the following: 

\vspace{0.1in}

\begin{theorem}~\label{thm:wreath}
Let $G \leq S_n$, and $G \wr S_2$ be the wreath product of $G$ with the symmetric group $S_2$.  Then
\[
i(P(G \wr S_2) ,t) = \sum_{k=0}^{t} i^2(P(G),k) \cdot i^2(P(G),t-k)
\]
for any integer $t \geq 2$.
\end{theorem}

\vspace{0.1in}

In Section~\ref{Miscellaneous}, we study miscellaneous permutation polytopes.  We begin by studying the permutation polytopes $P(A_n)$ where $A_n \leq S_n$, the \emph{alternating group} on $\{1,2,\ldots,n\}$, consists of permutations with even signature.  One of the main focuses in the literature is on determining the facets of $P(A_n)$.  Cunningham and Wang \cite{CunninghamWang2004}, and independently Hood and Perkinson \cite{HoodPerkinson2004}, proved that $P(A_n)$ has exponentially many facets in $n$, resolving a problem of Brualdi and Liu \cite{BrualdiLiu1991}.  However, a full facet description is still not known.  Moreover, no polynomial time algorithm in $n$ is known for membership in $P(A_n)$.  The difficulty of attaining a description of all facets of these polytopes is demonstrated by the following proposition, which shows that the first iterate of the theta body hierarchy for the polytopes $P(A_n)$ is almost never equal to $P(A_n)$ itself.

\vspace{0.1in}

\begin{proposition} \label{prp:non-two-level}
The polytope $P(A_n)$ is two-level, and hence $A_n$ is $TH_1$-exact, if and only if $n \leq 4$.  Moreover, for $n \geq 8$, $P(A_n)$ is at least $(\lfloor \frac{n}{4} \rfloor + 1)$-level.
\end{proposition}

\vspace{0.1in}

We conclude the paper with computations of volumes and Ehrhart polynomials of permutation polytopes for many subgroups of $S_3,S_4$, and $S_5$.

\section{Preliminaries}\label{Preliminaries}

Given a convex $d$-dimensional polytope $P \subset \R^n$, its \emph{Ehrhart polynomial}
$i(P,t)$ is the function that counts the number of points in
$\Z^n \cap tP$, where $tP = \{ tX \ | \ X \in P\}$ is the $t^{th}$
dilation of $P$. It is well
known that the volume of $P$ is the leading coefficient of the
polynomial $i(P,t)$.  In order to determine the Ehrhart polynomial or the normalized volume of
a lattice polytope $P$, it is often useful to know a unimodular
triangulation of $P$.  Let $L$ be an affine lattice $L \subset
\R^n$. A simplex with vertices $v_1, v_2, \ldots, v_m \in L$ is said to be
\emph{unimodular} in the lattice $L$ if $\{v_m-v_1, \ v_{m-1}-v_1,
\ldots, v_2-v_1\}$ is a basis for the lattice $L$.  A polytope whose
vertices lie in $L$ is said to have a \emph{unimodular triangulation}
in $L$ if it has a triangulation in which all maximal
dimensional simplices are unimodular in $L$.  We will be
interested in unimodularity with respect to two kinds of lattices.  We
say that a polytope $P$ is $P$-unimodular if it has a unimodular triangulation in
the lattice spanned by the vertices of $P$, and we will say that $P$ is $\Z$-unimodular if it has a unimodular triangulation
in the lattice aff$(P) \cap \Z^n$.  The following lemma shows that if
a polytope has a $\Z$-unimodular triangulation and
the number of faces of each dimension in the triangulation is known,
then its Ehrhart polynomial, and hence its volume, can be determined
immediately.

\begin{lemma} \label{lem:unimodularehrhart}(See Theorem 9.3.25 in \cite{DeLoeraRambauSantos2010})
Let $P \subset \R^n$ be a lattice polytope.  Assume that $P$ has a $\Z$-unimodular triangulation with $f_k$ faces of dimension $k$.  Then the Ehrhart polynomial of $P$ is 
\[
i(P,t) = \sum_{k=0}^n {t - 1 \choose k}f_k.
\]
\end{lemma}

In order to compute triangulations of our polytopes, we make repeated use of Gale duality.
In what follows, let $P \subset \R^n$ be a polytope with $r$ vertices $V=\{v_1,v_2,\ldots,v_r\}$ that lie on a common subspace and let $d = \dim(P)$.  Let $V \in \R^{n \times r}$ be the matrix given by
\begin{eqnarray}\label{affinization}
\begin{pmatrix}
v_1 & v_2 & \cdots & v_r
\end{pmatrix}.
\end{eqnarray}
Let $\mathcal{G} \in \R^{(r-d-1) \times r}$ be a matrix whose rows form a basis for the space of linear dependences of the columns of $(\ref{affinization})$.  The \emph{Gale dual} of $P$ is the vector configuration $\{\overline{v}_1,\overline{v}_2,\ldots,\overline{v}_r\}$ consisting of the columns of $\mathcal{G}$.  Note that $\mathcal{G}$ is unique up to linear coordinate transformations.  The relationship between triangulations of a polytope and the structure of its Gale dual hinges on the \emph{chamber complex} of $\mathcal{G}$.  Denote by $\Sigma_{\mathcal{G}}$ the set of cones generated by all bases of $\mathcal{G}$, that is, all subsets of  $\{\overline{v}_1,\overline{v}_2,\ldots,\overline{v}_r\}$ that form bases for the column space of $\mathcal{G}$.  If $\sigma \in \Sigma_{\mathcal{G}}$, let $\partial \sigma$ denote its boundary, and let $\partial \Sigma_{\mathcal{G}}$ be the union of the boundaries of all cones $\sigma \in \Sigma_{\mathcal{G}}$.  The complement of $\partial \Sigma_{\mathcal{G}}$ inside the cone generated by $\{\overline{v}_1,\overline{v}_2,\ldots,\overline{v}_r\}$ consists of open convex cones.  The closure of such an open convex cone is called a \emph{chamber}, and the \emph{chamber complex} of $\mathcal{G}$ is the collection of all these chambers.  The chamber complex of $\mathcal{G}$ and its relationship to triangulations of $P$ is encapsulated in the following lemma.

\begin{lemma}\label{lem:gale}(See Theorem 5.4.5, Theorem 5.4.7, and Theorem 5.4.9 in \cite{DeLoeraRambauSantos2010})
Let $P \subset \R^n$ be a $d$-dimensional polytope with vertex set $V = \{v_1,v_2,\ldots,v_r\}$ and Gale dual $\{\overline{v}_1, \overline{v}_2, \ldots, \overline{v}_r\}$.  Let $\tau$ be a chamber of the chamber complex of $\mathcal{G}$.  Then \[\Delta = \bigcup \mbox{conv}(V \backslash \{v_{j_1},v_{j_2},\ldots,v_{j_{r-d-1}}\} ),\] taken over all $\{v_{j_1},v_{j_2},\ldots,v_{j_{r-d-1}}\}$ such that $\tau \subseteq \mbox{conv}\{\overline{v}_{j_1},\overline{v}_{j_2},\ldots,\overline{v}_{j_{r-d-1}}\}$ is a full-dimensional cone in the Gale dual, is a regular triangulation of $P$.  Moreover, all regular triangulations of $P$ arise in this way from some chamber $\tau$.
\end{lemma}

Unfortunately, the aforementioned triangulations given by the Gale dual may not be $\Z$- nor $P$-unimodular, so we still need methods to determine if a given polytope $P$ has a $\Z$-unimodular or $P$-unimodular triangulation. One way to do this is through the use of Gr\"{o}bner bases of toric ideals.  Though this can be addressed in a more general setting, we will restrict ourselves to permutation polytopes arising from subgroups of a particular $S_n$.  Let $G = \{g_1,g_2,\ldots,g_k\}$ be elements of such a subgroup, and as usual consider $g_i$ as an $n \times n$ permutation matrix for each $i$.  Let $\C[\textbf{x}] = \C[x_{g_1},x_{g_2},\ldots,x_{g_k}]$ be the polynomial ring in $k$ indeterminates indexed by the elements of $G$ and let $\C[\textbf{t}] := \C[t_{\ell m} : 1 \leq \ell,m \leq n]$.  The algebra homomorphism induced by the map
\[
\hat{\pi}_G : \C[\textbf{x}] \to \C[\textbf{t}], \ \hat{\pi}_G(x_{g_i}) = \prod_{1 \leq \ell,m \leq n } t_{\ell m}^{(g_i)_{\ell m}}, \ \ 1 \leq i \leq k
\]  
has as its kernel the ideal $I_G$.  Given a monomial order $\prec$ on $\C[\textbf{x}]$, the ideal $I_G$ can determine a $P(G)$-unimodular triangulation of $P(G)$.  Moreover, this triangulation is always regular.  See \cite{DeLoeraRambauSantos2010,Sturmfels1996} for more on regular triangulations.

\begin{lemma}\label{lem:unimodulargrobner}(See Corollary 8.9 in \cite{Sturmfels1996} and Theorem 9.4.5 in \cite{DeLoeraRambauSantos2010})
Let $in_{\prec}(I_G)$ be the initial ideal of $I_G$ with respect to the term order $\prec$.  The support vectors of the generators of the radical of $in_{\prec}(I_G)$ are the minimal non-faces of a regular triangulation of $P(G)$.  Moreover, $in_{\prec}(I_G)$ is square-free if and only if the corresponding triangulation $\Delta_{\prec}$ of $P(G)$ is $P(G)$-unimodular.
\end{lemma}
\noindent By Lemma~\ref{lem:unimodulargrobner} and the theory of Gr\"{o}bner bases, $P(G)$ will have a $P(G)$-unimodular triangulation if there is a term order $\prec$ on $\C[\textbf{x}]$ such that the Gr\"{o}bner basis of $I_G$ is generated by polynomials whose initial terms are square-free.  This will be exploited in Section~\ref{Frobenius}.  For more on the relationship between toric ideals, Gr\"{o}bner bases, and triangulations, see \cite{Sturmfels1996}.

Recently, Gouveia, Parrilo, and Thomas \cite{GouveiaParriloThomas2010} constructed a hierarchy of convex bodies, each given as the projection of the feasible region of a semidefinite program, approximating the convex hull of an arbitrary real variety.  In Section~\ref{Frobenius}, given a permutation polytope $P(G) \subset \R^{n \times n}$, we discuss the convergence of this hierarchy of relaxations $TH_1(I) \supseteq TH_2(I) \supseteq \cdots \supseteq P(G)$, which are known as theta bodies. Here, $I$ is an ideal in $\R[x_{ij} : 1 \leq i, j \leq n]$ whose real variety is the vertex set of $P(G)$.  This hierarchy of relaxations has the property that if $TH_k(I) = P(G)$ for some fixed $k$, linear optimization over $P(G)$ can be performed in polynomial time in $n$ provided a certain algebraic oracle.  Of particular interest then are polytopes for which $TH_1(I) = P(G)$, in which case we say $G$ is \emph{$TH_1$-exact}.  The concept of $TH_1$-exactness was defined in \cite{GouveiaParriloThomas2010} for general ideals $I$, but we focus on ideals whose zero-sets are vertices of permutation polytopes.  $TH_1$-exact varieties are characterized polyhedrally in \cite{GouveiaParriloThomas2010}, and we again restrict this characterization to permutation polytopes.

\begin{lemma}\label{lem:1exact}(See Theorem 4.2 in \cite{GouveiaParriloThomas2010})
The group $G$ is $TH_1$-exact if and only if for every facet defining inequality $c \cdot x - \alpha \geq 0$ of $P(G)$, there is a plane $c \cdot x - \beta = 0$ parallel to $c \cdot x - \alpha =0$ such that all vertices of $P(G)$ lie in $\{x \ | \ c \cdot x - \alpha = 0\} \cup \{x \ | \ c \cdot x - \beta = 0\}$. 
\end{lemma}

A polytope satisfying the facet property in Lemma~\ref{lem:1exact} is called \emph{compressed} or \emph{two-level}. Using Lemma~\ref{lem:1exact}, we can use Gale duality to characterize groups whose permutation polytopes are $TH_1$-exact.  We do so in the following lemma, which was proved independently by Gouveia, Parrilo, and Thomas, but we provide a self-contained proof.

\begin{lemma}\label{lem:galeexact}
Let $P(G) \subseteq \R^{n \times n}$ be a permutation polytope with vertex set $\{v_1,v_2,\ldots,v_r\}$ and Gale dual $\{\overline{v}_1,\overline{v}_2,\ldots,\overline{v}_r\}$.  Then $G$ is $TH_1$-exact if and only if for every $J \subseteq \{1,2,\ldots,r\}$ such that $\mbox{conv} \{ v_j \ | \ j \in J \} $ is a facet of $P(G)$, $\sum_{j \notin J} \overline{v}_j = 0$.
\end{lemma}

\proof
Throughout this proof, we use the equivalence of $TH_1$-exactness with the property of $P(G)$ being two-level, which was proved in Lemma~\ref{lem:1exact}.  Let $J \subseteq \{1,2,\ldots,r\}$ such that conv$\{v_j \ | \ j \in J\}$ is a facet of $P(G)$ with the defining inequality $c \cdot x - \alpha \geq 0$ valid on $P(G)$.  Then
\[
0 = (c, -\alpha) \begin{pmatrix} v_1 & v_2 & \cdots & v_r \\ 1 & 1 & \cdots & 1 \end{pmatrix} \begin{pmatrix} \overline{v}_1 \\ \overline{v}_2 \\ \vdots \\ \overline{v}_r \end{pmatrix} = (c \cdot v_1 -\alpha, c \cdot v_2 - \alpha, \cdots, c \cdot v_r - \alpha) \begin{pmatrix} \overline{v}_1 \\ \overline{v}_2 \\ \vdots \\ \overline{v}_r \end{pmatrix}.
\]
Since $G$ is exact, $c \cdot v_j - \alpha$ can take at most two values.  By construction, one of these values is zero. If the other value is $\beta$, then $c \cdot v_j - \alpha = \beta$ if and only if $j \notin J$.  Thus $\sum_{j \notin J} \beta \overline{v}_j = 0$, which implies that $\sum_{j \notin J} \overline{v}_j = 0$ since $\beta \neq 0$.

For the converse, suppose that $\sum_{j \notin J} \overline{v}_j = 0$ for every $J$ such that conv$\{v_j \ | \ j \in J\}$ is a facet.  Fix such a $J$ and assume that the facet inequality of $P(G)$ defining it is $c \cdot x - \alpha \geq 0$.  Then, as we have done above,
\[
0 = \sum_{j \notin J} (c \cdot v_j - \alpha) \overline{v}_j = \sum_{j \notin J} (c \cdot v_j) \overline{v}_j - \alpha \sum_{j \notin J} \overline{v}_j = \sum_{j \notin J} (c \cdot v_j) \overline{v}_j.
\]
Suppose that there are at least two distinct values among $\{c \cdot v_j \ | \ j \notin J\}$ and that $\gamma$ is the least value.  Then
\[
0 = \sum_{j \notin J} (c \cdot v_j -\gamma) \overline{v}_j
\]
yields a positive dependence relation on $\{\overline{v}_j \ | \ j \notin J\}$ that does not use all the elements in the set.  This contradicts the assumption that $J$ induces a facet of $P(G)$.  Thus $\{c \cdot v_j \ | \ j \notin J\}$ has only one element, and hence $G$ is exact.
\qed

Note that in particular if $P(G)$ is $TH_1$-exact, then it contains a $P(G)$-unimodular triangulation.  Moreover, as we will see, all simplices in this triangulation have the same volume.

\section{Cyclic and Dihedral Groups}~\label{cyclicdihedral}

We dedicate this section to the proof of Theorem~\ref{thm:cyclicdihedralehrhart}. We begin with the following lemma, which determines the Ehrhart polynomial of $P(C_n)$ and proves part (1) of Theorem~\ref{thm:cyclicdihedralehrhart}.

\begin{lemma}\label{lem:cyclic}
The Ehrhart polynomial of $P(C_n)$ is $i(P(C_n),t) = \binom{t + n-1}{n-1}$.
\end{lemma}

\proof
Let $t$ be a positive integer and let $\phi: t P(C_n) \to t \Delta_n$ be the affine map given by $\phi(X) = [X_{1,1}, \ X_{1,2}, \ \cdots \ X_{1,n}]^T$.  Here, $\Delta_n$ is the standard $(n-1)$-simplex $\mbox{conv}\{e_i \ | \ 1 \leq i \leq n\} \subseteq \R^n$.  Note that $\phi$ is a well-defined map since the sum of the first row of any matrix in $tP(C_n)$ is $t$.  We claim that $\phi$ induces a bijection between the sets $tP(C_n) \cap \Z^{n \times n}$ and $t \Delta_n \cap \Z^{n \times n}$.  If $X \in tP(C_n)$ is an integer matrix, then its first row contains integer entries whose sum is $t$,  so $\phi(X)$ is indeed an integer point in $t \Delta_n$.  It suffices to show that every integer point in $t \Delta_n$ has a unique integral pre-image.  Let $[X_{1,1}, \ X_{1,2}, \cdots X_{1,n}]^T \in t \Delta_n \cap \Z^n$.  Then $\phi( X_{1,1} \cdot e + X_{1,2} \cdot  g + \cdots + X_{1,n} \cdot g^{n-1}) =  [X_{1,1}, \ X_{1,2}, \cdots X_{1,n}]^T$, so $[X_{1,1}, \ X_{1,2}, \cdots X_{1,n}]^T$ has an integral pre-image.  Moreover, this pre-image is unique, since $g^i$ is the only vertex of $tP(C_n)$ whose $(1,i+1)$-entry is non-zero.  Thus $\phi$ induces a bijection between $tP(C_n) \cap \Z^{n \times n}$ and $t \Delta_n \cap \Z^{n \times n}$, and the result follows since $i(\Delta_n,t) = \binom{t+n-1}{n-1}$.  The fact that the volume of $P(C_n)$ is $\frac{1}{(n-1)!}$ follows because the first coefficient of $i(P(C_n),t)$ is $\frac{1}{(n-1)!}$ and $\Delta_n$ is $\Z$-unimodular. 
\qed

We now investigate the polytopes $P(D_n)$ and their Ehrhart polynomials. Recall the following lemma concerning the dimension of $P(D_n)$.

\begin{lemma}\label{lem:dihedraldim}(See Theorem 4.1 of \cite{Steinkamp1999})
The dimension of the polytope $P(D_n)$ is $2n-2$ if $n$ is odd and $2n-3$ if $n$ is even.
\end{lemma}

Lemma~\ref{lem:dihedraldim} indicates that Gale duality is very useful for determining the Ehrhart polynomial of $P(D_n)$, since the Gale dual lies in a space of dimension $|D_n| - \dim (P(D_n))-1$, which is one if $n$ is odd and two if $n$ is even.

\begin{lemma}\label{lem:dihedralgale}
If $n$ is odd, the Gale dual of $P(D_n)$ is a vector configuration in $\R$ consisting of $n$ copies of each of the vectors $\pm 1$.  If $n$ is even, $n=2m$, the Gale dual of $P(D_n)$ is the vector configuration in $\R^2$ consisting of $m$ copies of each of the four vectors $[\pm 1, 0]^T,[0,\pm 1]^T$.
\end{lemma}

\proof
Throughout this proof, let $\mathcal{G}$ be the matrix whose columns form the Gale dual of $P(D_n)$ with its columns indexed by $\{e,r,r^2,\ldots,r^{n-1},f,fr,fr^2,\ldots,fr^{n-1}\}$ in that order.  The following linear relation holds for $D_n$:
\begin{eqnarray}\label{galesum}
e + r + r^2 + \ldots + r^{n-1} = f + fr + fr^2 + \ldots + fr^{n-1} = J_{n \times n},
\end{eqnarray}
where $J_{n \times n}$ is the $n \times n$ matrix whose entries are all one.  When $n$ is odd, Lemma~\ref{lem:dihedraldim} implies that $P(D_n)$ is $2n-2$ dimensional, so the Gale dual of $P(D_n)$ is one dimensional.   Thus, Equation~(\ref{galesum}) implies that
\[
\mathcal{G} = \begin{pmatrix}
1 \ 1 \ \cdots \ 1 \ -1 \ -1 \ \cdots \ -1
\end{pmatrix},
\]
with $n$ copies of $1$ and $n$ copies of $-1$.  When $n$ is even, $n=2m$, Lemma~\ref{lem:dihedraldim} implies that $P(D_n)$ is $2n-3$ dimensional, so the Gale dual of $P(D_n)$ is two dimensional.
We observe that the relation
\begin{eqnarray}\label{galesum2}
\sum_{j=0}^{m-1} r^{2j+1} = \sum_{j=0}^{m-1} f r^{2j}
\end{eqnarray}
holds for $D_n$ when $n$ is even. The linear relations (\ref{galesum})-(\ref{galesum2}) and (\ref{galesum2}) are linearly independent, so we deduce that
\[
\mathcal{G} = \begin{pmatrix}
1 & 0 & \ldots & 0 & 0 & -1 & \ldots & -1  \\
0 & 1 & \ldots & 1 & -1 & 0 & \ldots & 0 
\end{pmatrix}.
\]
We conclude that the Gale dual is the vector configuration in $\R^2$ consisting of $n$ copies of each of the four vectors $[\pm 1,0]^T, [0, \pm1]^T$.
\qed

We now compute the Ehrhart polynomial of $P(D_n)$.  The symmetry in its Gale dual shows that the number of faces of a given dimension in any regular triangulation of $P(D_n)$ is the same.  Note that this in principle completely describes the secondary polytope of $P(D_n)$.  Thus, if we find any $\Z$-unimodular triangulation of $P(D_n)$ and compute the number of faces of each dimension in any other regular triangulation of it, we can recover its Ehrhart polynomial via Lemma~\ref{lem:unimodularehrhart}.  We begin by finding a $P(D_n)$-unimodular triangulation.

\begin{proposition}\label{prp:dihedralunimodular}
The polytope $P(D_n)$ has a $P(D_n)$-unimodular regular triangulation.
\end{proposition}

\proof
Let $G$ be the graph with vertices  $\{1,2,\ldots,n\}$ and edges $i,i+1$ for each $i \in \{1,2,\ldots,n\}$. Let $A_G$ be the adjacency matrix of $G$. Consider the polytope
\[
P_G = \left \{ X \in [0,1]^{n \times n} \ : \ A_G X = X A_G, \  \ \sum_{j=1}^n X_{ij} = 1 \ \forall i, \ \ \sum_{i=1}^n X_{ij} = 1 \ \forall j \right \}.
\]
The integer points of $P_G$ are permutations commuting with $A_G$, so they are precisely the automorphisms of $G$.  Since the automorphism group of $G$ is $D_n$ and $P_G$ is integral (see Theorem 2 of \cite{Tinhofer1986}), this implies that $P_G = P(D_n)$.  But by Theorem 4.4 of \cite{DeLoeraHillarMalkinOmar2010}, the vertex set of $P_G$ is exact, which by Theorem 2.4 of \cite{Sullivant2006} and Theorem 4.2 of \cite{GouveiaParriloThomas2010} implies that every reverse lexicographic triangulation of $P(D_n)$ is $P(D_n)$-unimodular. Since reverse lexicographic triangulations are regular, the result follows.
\qed

\vspace{0.1in}

In order to establish that $P(D_n)$ is $\Z$-unimodular, we prove that the index of the lattice generated by its vertices in the lattice aff$(P(D_n)) \cap \Z^{n \times n}$ is one.
\begin{proposition}\label{prp:dihedralindex}
The index of the lattice generated by the vertices of $P(D_n)$ has index one in the lattice $\mbox{aff}(D_n) \cap \Z^{n \times n}$.
\end{proposition}
\proof
First, consider when $n$ is odd.  For simplicity, let $D_n$ consist of the matrices $v_1,v_2,\ldots,v_{2n}$, where $v_{2i+1} = r^i$ and $v_{2i+2}$ is the unique flip in $D_n$ fixing $i+1$, $0 \leq i \leq n-1$.  It suffices to prove that if $X \in \Z^{n \times n}$ and $X$ is an $\R$-linear combination of the matrices $\{v_{2n} - v_1,v_{2n-1}-v_1,\ldots,v_2-v_1\}$, then $X$ is a $\Z$-linear combination of these matrices.  Assume then that $X = \sum_{j=2}^{2n} \alpha_j (v_1 - v_j) = \left( \sum_{j=2}^{2n} \alpha_j \right) v_1- \sum_{j=2}^{2n} \alpha_j v_j$, with $\alpha_j \in \R$.  Let $\alpha \in [0,1)$ such that $\sum_{j=2}^{2n} \alpha_j + \alpha \in \Z$.  Since $e$ and $v_{2i+2}$ are the only elements of $D_n$ with the $(i+1,i+1)$-entry in their support, and since $X$ has integer entries, we conclude that $\alpha_{2i+2} - \alpha \in \Z$ for all $0 \leq i \leq n-1$.  Moreover, for any $i$, there is a unique flip with the $(1,i+1)$-entry in its support.  Since $r^i$ is the only rotation with the $(1,i+1)$-entry in its support, and again since $X$ has integer entries, we deduce that $\alpha_{2i+1} + \alpha \in \Z$ for all $0 \leq i \leq n-1$.  Now recall from Equation~(\ref{galesum}) in the proof of Lemma~\ref{lem:dihedralgale} that $\sum_{i=0}^{n-1} v_{2i+1} - \sum_{i=0}^{n-1} v_{2i+2} = 0$, so we have that
\begin{align*}
X & = \left( \sum_{j=2}^{2n} \alpha_j \right) v_1- \sum_{j=2}^{2n} \alpha_j v_j - \alpha \left( \sum_{i=0}^{n-1} v_{2i+1} - \sum_{i=0}^{n-1} v_{2i+2} \right)  \\
&= \sum_{i=1}^{n-1} (\alpha_{2i+1} + \alpha) (v_1 - v_{2i+1}) + \sum_{i=0}^{n-1} (\alpha_{2i+2} - \alpha) (v_1 - v_{2i+2}),
\end{align*}
and hence $X$ is a $\Z$-linear combination of $\{v_{2n} - v_1,v_{2n-1}-v_1,\ldots,v_2-v_1\}$.

Now consider when $n$ is even, $n=2m$.  We let $D_n$ consist of the $4m$ vectors $\{u_1,u_2,\ldots,u_{2m}\}$, $\{v_1,v_2,\ldots,v_{2m}\}$, where $u_{2i+1} = r^{2i}$, $v_{2i+1} =  r^{2i+1}$, $u_{2i+2}$ is the unique flip supported on the $(1,2i+1)$-entry, and $v_{2i+1}$ is the unique flip supported on the $(1,n+2i)$-entry, entries taken mod $i$ and $0 \leq i \leq m-1$.  Suppose that $X$ is an $\R$-linear combination of the form
\[
X = \sum_{i=2}^{2m} \alpha_i (u_1 - u_i) + \sum_{i=1}^{2m} \beta_i (u_1 - v_i) = \left(\sum_{i=2}^{2m} \alpha_i + \sum_{i=1}^{2m} \beta_i\right)u_1 - \sum_{i=2}^{2m} \alpha_i u_i - \sum_{i=1}^{2m} \beta_i v_i.
\]
Let $\alpha \in [0,1)$ such that $\sum_{i=2}^{2m} \alpha_i + \sum_{i=1}^{2m} \beta_i - \alpha \in \Z$.  Notice that $u_1=e$.  Since $e$ and $u_{2i+2}$ are the only elements of $D_n$ with the $(i,i)$-entry in their support and $X$ has integer entries, we conclude that $\alpha_{2i_2} - \alpha \in \Z$ for all $0 \leq i \leq m-1$.  Moreover, $u_{2i+1}$ and $u_{2i+2}$ are the only elements with the $(1,2i+1)$-entry in their support, so $\alpha_{2i+1} + \alpha \in \Z$ for all $0 \leq i \leq m-1$.  Similarly, if $\beta \in [0,1)$ such that $\beta_1 + \beta \in \Z$, then $\beta_{2i+1} + \beta, \beta_{2i+2} - \beta \in \Z$ for all $0 \leq i \leq m-1$.  Now from Equation (\ref{galesum})-Equation (\ref{galesum2}) and Equation (\ref{galesum2}) in the proof of Lemma~\ref{lem:dihedralgale}, we have that
\[
\sum_{i=0}^{m-1} u_{2i+1} = \sum_{i=0}^{m-1} u_{2i+2}, \ \  \sum_{i=0}^{m-1} v_{2i+1} = \sum_{i=0}^{m-1} v_{2i+2}.
\] 
We conclude then that 
\begin{align*}
X &= \left(\sum_{i=2}^{2m} \alpha_i + \sum_{i=1}^{2m} \beta_i \right) u_1 - \sum_{i=2}^{2m} \alpha_i u_i - \sum_{i=1}^{2m} \beta_i v_i - \alpha \left( \sum_{i=0}^{m-1} u_{2i+1} - \sum_{i=0}^{m-1} u_{2i+2} \right) \\ & - \beta \left( \sum_{i=0}^{m-1} v_{2i+1} - \sum_{i=0}^{m-1} v_{2i+2} \right) \\
&=  \sum_{i=1}^{m-1} (\alpha_{2i+1} + \alpha) (u_1 - u_{2i+1}) + \sum_{i=0}^{m-1} (\alpha_{2i+2} - \alpha) (u_1 - u_{2i+2}) +  \sum_{i=0}^{m-1} (\beta_{2i+1} + \beta) (u_1 - v_{2i+1}) \\ &+ \sum_{i=0}^{m-1} (\beta_{2i+2} - \beta) (u_1 - v_{2i+2}),
\end{align*}
which is a $\Z$-linear combination of the required vectors.
\qed


We now determine the number of faces of each dimension in a particular triangulation of $P(D_n)$.  This together with Lemma~\ref{lem:unimodularehrhart}, Proposition~\ref{prp:dihedralunimodular}, and Proposition~\ref{prp:dihedralindex} proves parts (2) and (3) of Theorem~\ref{thm:cyclicdihedralehrhart}.

\vspace{0.1in}

\proof [Theorem \ref{thm:cyclicdihedralehrhart}]
First consider when $n$ is odd.  By Lemma~\ref{lem:dihedralgale}, the Gale dual of $P(D_n)$ consists of the vectors $\{e_1^{(1)},e_1^{(2)},\ldots,e_1^{(n)},-e_1^{(1)},-e_1^{(2)},\ldots,-e_1^{(n)}\}$ where the $e_1^{(i)},-e_1^{(i)}$ are copies of the vectors $e_1,-e_1$ in $\R$ respectively, $1 \leq i \leq n$. The set consisting of the vector $e_1$ is the only extreme ray in one chamber in the Gale dual, so by Lemma~\ref{lem:gale}, $P(D_n)$ has a triangulation $\Delta$ with maximal dimensional simplices $\left \{ \mbox{conv} \{ G \backslash \{r^i\} \}  | \ 1 \leq i \leq n \right \}$.  The number of $(k+1)$-element subsets of $G$ is $\binom{2n}{k+1}$.  By Lemma~\ref{lem:gale}, of these subsets, the ones that are not simplices in $\Delta$ are those that contain all of $\{e,r,r^2,\ldots,r^{n-1}\}$.  There are precisely $\binom{2n-n}{k+1-n}$ such subsets, so we conclude that the number of $k$-dimensional faces $f_k$ in $\Delta$ is
\[
f_k = \binom{2n}{k+1} - \binom{n}{k+1-n}.
\]
By the symmetry in the Gale dual, this is also the number of $k$-dimensional faces in a reverse lexicographic triangulation of $P(D_n)$, which is $P(D_n)$-unimodular by Proposition~\ref{prp:dihedralunimodular} and hence $\Z$-unimodular by Propostion~\ref{prp:dihedralindex}. The Ehrhart polynomial follows from Lemma~\ref{lem:unimodularehrhart}.  Moreover, we see that $f_{2n-2} = \binom{2n}{2n-1} - \binom{n}{n-1} = n$, so the volume of $P(D_n)$ is $\frac{n}{(2n-2)!}$.

\vspace{0.1in}

Now consider $P(D_n)$ when $n$ is even, $n=2m$.  By Lemma~\ref{lem:dihedralgale}, the Gale dual of $P(D_n)$ consists of the copies $\{e_1^{(i)}, e_2^{(i)}, -e_1^{(i)}, -e_2^{(i)} \ | \ 1 \leq i \leq m\}$ of $e_1,e_2,-e_1,-e_2$ respectively in $\R^2$.  Consider the chamber of the Gale dual whose extreme rays are the vectors $\{e_1,e_2\}$.  By Lemma~\ref{lem:gale}, this chamber gives the regular triangulation $\Delta$ of $P(D_n)$ whose maximal dimensional simplices are $\left \{ \mbox{conv} \{ G \backslash \{ r^{2i-1}, r^{2j} \} \}  | \ 1 \leq i,j \leq n \right \}$.  By a similar counting argument as in the odd case, we conclude that
\[
f_k = \binom{2n}{k+1} - \binom{2}{1} \binom{2n-m}{k+1-m} + \binom{2n-2m}{k+1-2m}.
\] 
Again, since Lemma~\ref{lem:unimodularehrhart} implies that $P(D_n)$ has a $P(D_n)$-unimodular triangulation and hence by Proposition~\ref{prp:dihedralunimodular} a $\Z$-unimodular triangulation with the same face numbers, the Ehrhart polynomial follows by Lemma~\ref{lem:unimodularehrhart}.  Lastly, we see that the volume of $P(D_n)$ when $n$ is even is $f_{2n-3}$, which is
\[
\frac{1}{(2n-3)!} \left( \binom{2n}{2n-2} - 2 \binom{2n-m}{2n-m-2} + \binom{2n-m}{2n-2m-2} \right) = \dfrac{n^2}{4 \cdot (2n-3)!}.
\]
\qed

\noindent \textbf{Remark.} One can further show from the proof of Theorem~\ref{thm:cyclicdihedralehrhart} that all simplices in the given triangulations have the same volume.

\vspace{0.1in}

\begin{example}
We use Theorem~\ref{thm:cyclicdihedralehrhart} to determine the Ehrhart polynomials and volumes of $P(D_4)$ and $P(D_5)$.  First, we have that 
$i({P(D_4)},t)$ is
\[
8 \binom{t-1}{0} + 26 \binom{t-1}{1} + 44 \binom{t-1}{2} + 41 \binom{t-1}{3} + 20 \binom{t-1}{4} + 4 \binom{t-1}{5},
\]
which is precisely the polynomial
\[
\frac{1}{30}t^5 + \frac{1}{3}t^4 + \frac{4}{3}t^3 + \frac{8}{3} t^2 + \frac{79}{30}t + 1.
\]
The volume of $P(D_4)$ is therefore $\frac{1}{30}$.  Similarly, we have that $i({P(D_5)},t)$ is
\begin{align*}
& 10 \binom{t-1}{0} + 45 \binom{t-1}{1} + 120 \binom{t-1}{2} + 210 \binom{t-1}{3} + 251 \binom{t-1}{4} + \\
& 205 \binom{t-1}{5} + 110 \binom{t-1}{6} + 35 \binom{t-1}{7} + 5 \binom{t-1}{8},       
\end{align*}
which is precisely the polynomial
\[
\frac{1}{8064}t^8 + \frac{5}{2016}t^7 + \frac{5}{192}t^6 + \frac{25}{144}t^5 + \frac{95}{128}t^4 + \frac{575}{288}t^3 + \frac{6515}{2016}t^2 + \frac{475}{168}t + 1.
\]
The volume of $P(D_5)$ is therefore $\frac{1}{8064}$.
\end{example}

\section{Frobenius Groups}~\label{Frobenius}

In this section, we discuss triangulations and normalized volumes of Frobenius polytopes, leading to a proof of Theorem~\ref{thm:Frobenius}.  We also establish that all Frobenius groups are exact, hence proving Proposition~\ref{prp:frobeniustwolevel}.  For the remainder of this section, we assume that $G \leq S_n$ is a Frobenius group.  We let $N = \{u_1,u_2,\ldots,u_{n}\}$ be its Frobenius kernel ($n = |N|$), and we let $H=\{v_1,v_2,\ldots,v_{h}\}$ be its Frobenius complement ($h = |H|$).  We assume throughout that $H$ is the set of coset representatives for $N$ in $G$.  We let $\mathcal{G}$ denote the matrix whose columns form the Gale dual of $P(G)$.  Recall that $G = NH$ and $H \cap N = \{e\}$, and so $G$ consists of the $nh$ matrices
\[
u_1v_1, u_2v_1, \ldots, u_n v_1, u_1v_2, u_2v_2, \ldots, u_nv_2, \ldots, u_1v_h, u_2v_h, \ldots, u_n v_h
\]
and we index the columns of $\mathcal{G}$ by $G$ in this order.  The following lemmas are proven in \cite{CollinsPerkinson2008}.

\begin{lemma}\label{lem:frobeniusequations}(See Proposition 4.2 in \cite{CollinsPerkinson2008})
If $G \leq S_n$ is Frobenius, then $\sum_{i=1}^{n} u_i v_j = J_{n \times n}$ for all $j$, $1 \leq j \leq h$, where $J_{n \times n}$ is the $n \times n$ matrix of all 1s.
\end{lemma}

\begin{lemma}\label{lem:frobeniusdimension}(See Corollary 4.5 in \cite{CollinsPerkinson2008})
If $G \leq S_n$ is Frobenius, the dimension of $P(G)$ is $|G|-|H|$.
\end{lemma}

Lemma~\ref{lem:frobeniusequations} gives us the $|H|-1$ linearly independent relations $\sum_{i=1}^n u_i v_1 = \sum_{i=1}^n u_i v_j$, $2 \leq j \leq h$.  The dimension formula in Lemma~\ref{lem:frobeniusdimension} tells us that the $|H|-1$ relations in Lemma~\ref{lem:frobeniusequations} actually form a basis for the space of linear dependences of $G$.  As a consequence, we get the Gale dual of $P(G)$.

\begin{proposition}\label{prop:frobeniusgale}
The Gale dual of $P(G)$ is the vector configuration consisting of $n$ copies $\{{\bf{1}}^{(1)},{\bf{1}}^{(2)},\ldots,{\bf{1}}^{(n)}\}$ of the all-ones vector $\bf{1}$ in $\R^{h-1}$, together with $n$ copies $\{-e_i^{(1)},-e_i^{(2)},\ldots,-e_i^{(n)}\}$ in $\R^{h-1}$ of $-e_i$ for $1 \leq i \leq h-1$, where $e_i$ is the $i^{th}$ standard basis vector.  In particular, the $u_iv_j$ column of the matrix $\mathcal{G}$ is the vector $\bf{1}$ if $j=1$, and $-e_{j-1}$ otherwise.
\end{proposition}

\proof
This follows directly from Lemma~\ref{lem:frobeniusequations} and Lemma~\ref{lem:frobeniusdimension}.
\qed

\vspace{0.1in}

Now consider the chamber in the Gale dual whose extreme rays are $\{-e_1,-e_2,\ldots,-e_{h-1}\}$.  From Lemma~\ref{lem:gale}, $P(G)$ has a corresponding regular triangulation $\Delta$ whose maximal dimensional simplices are
\begin{eqnarray}\label{eqn:triangulation}
\Delta = \left \{ \mbox{conv} \left \{ G \backslash \{u_{i_1} v_2, u_{i_2}v_3, \ldots, u_{i_{h-1}}v_h\} \right \} \ | \ 1 \leq i_j \leq n \right \}
\end{eqnarray}

Furthermore, from the structure of the Gale dual as given by Proposition~\ref{prop:frobeniusgale}, all triangulations of $P(G)$ have the same number of $k$-dimensional faces for any $k$.  Thus, if we can determine a $P(G)$-unimodular triangulation of $P(G)$ and count the number of faces of dimension $k$ for each $k$ in the triangulation $\Delta$, we can prove Theorem~\ref{thm:Frobenius}.  We proceed by showing that $P(G)$ has a $P(G)$-unimodular triangulation and then by determining the number of faces of given dimensions in $\Delta$.

\vspace{0.1in}

\begin{proposition}\label{prop:frobeniusunimodular}
If $G$ is Frobenius then $G$ has a $P(G)$-unimodular triangulation.
\end{proposition}

\proof
Our proof appeals to toric algebra.  Let $A \in \R^{n^2 \times |G|}$ be the matrix whose columns are the elements of $G$ written as $n^2$-dimensional column vectors by reading rows left to right and top to bottom.  We index the columns of $A$ by the elements of $G$ as in $\mathcal{G}$.  The toric ideal $I_G \subseteq \C[\textbf{x}] = \C[x_{u_rv_s}  : 1 \leq r \leq n, \ 1 \leq s \leq h]$ is the kernel of the homomorphism
\[
\hat{\pi}: \C[\textbf{x}] \to \C[\textbf{t}], \ \ \hat{\pi}(x_{u_rv_s}) = \prod_{1 \leq \ell,m \leq n} t_{\ell m}^{{(u_rv_s)}_{\ell m}}, 
\]
and by Lemma 4.1 of \cite{Sturmfels1996}, $I_G = \langle x^u - x^v \ | \ A(u-v) = 0, \ u,v \in \Z^{|G|} \rangle$.  By Lemma~\ref{lem:frobeniusequations} and Lemma~\ref{lem:frobeniusdimension}, $\ker(A)$ has the basis $\{b_1,b_2,\ldots,b_k\}$ where $b_i = e_{u_1} + e_{u_2} + \cdots + e_{u_n} - e_{u_1v_i} - e_{u_2v_i} - \cdots - e_{u_nv_i}$ for each $i$. Now if $u-v \in \ker(A)$ is integral, then $u-v = \sum_{i=1}^h \lambda_i b_i$, where $\lambda_i \in \Q$ for each $i$.  In fact, $\lambda_i \in \Z$ for each $i$ since the $u_{\ell}v_i$ component of $u-v$ is $\pm \lambda_i$.  We conclude by Corollary 4.4 of \cite{Sturmfels1996} that $I_G = \langle x_{H_1} - x_{H_{\ell}} : 2 \leq \ell \leq h \rangle$ where $x_{H_{\ell}} = \prod_{i=1}^{n} x_{u_i v_{\ell}}$ for each $\ell$.  

In fact, $\{x_{H_1} - x_{H_{\ell}} : 2 \leq \ell \leq h \}$ is a Gr\"{o}bner basis for $I_G$ with respect to the reverse lexicographic order $\prec$; here, $u_{r_1}v_{s_1}$ comes lexicographically before $u_{r_2}v_{s_2}$ if and only if $r_1 \leq r_2, s_1 \leq s_2$.  To see this, we use Buchberger's algorithm. For an introduction to this algorithm and details of terms to follow, see \cite{CoxLittleO'Shea2007}.  Consider any pair of polynomials $f_r = x_{H_1} - x_{H_r}, f_s = x_{H_1} - x_{H_s}$ in our generating set for $I_G$.   With respect to $\prec$, we compute the $S$-pair $S(f_r,f_s)$ and see that
\[
S(f_r,f_s) = \frac{x_{H_r}x_{H_s}}{-x_{H_r}} (x_{H_1} - x_{H_r}) - \frac{x_{H_r}x_{H_s}}{-x_{H_s}} (x_{H_1} - x_{H_s}) = x_{H_1}x_{H_r} - x_{H_1}x_{H_s}.
\]
Now since $x_{H_1}x_{H_r} - x_{H_1}x_{H_s} = x_{H_1}(x_{H_1} - x_{H_s}) - x_{H_1}(x_{H_1} - x_{H_r})$, we see that $\overline{S(f_r,f_s)}^{f_rf_s}=0$.  Since $r,s$ were arbitrary, Buchberger's algorithm concludes that $\{x_{H_1} - x_{H_{\ell}} : 2 \leq \ell \leq h \}$ is a Gr\"{o}bner basis for $I_G$.  By Lemma~\ref{lem:unimodulargrobner}, we conclude that $P(G)$ has a $P(G)$-unimodular triangulation.
\qed

We now proceed to prove Theorem~\ref{thm:Frobenius}.

\vspace{0.1in}

\proof [Theorem~\ref{thm:Frobenius}]
By Proposition~\ref{prop:frobeniusunimodular}, $P(G)$ has a $P(G)$-unimodular triangulation, and by Proposition~\ref{prop:frobeniusgale} and Lemma~\ref{lem:gale}, all triangulations of $P(G)$ have the same face numbers.  Thus it suffices to determine the number of top dimensional faces in the triangulation $\Delta$ in (\ref{eqn:triangulation}) and apply Lemma~\ref{lem:unimodularehrhart}.  We more generally determine the number of $k$-dimensional faces $f_k$ for each $k$.  Any $k$-simplex in $\Delta$ must be a subset of some maximal dimensional simplex of $\Delta$, and by Lemma~\ref{lem:gale}, all maximal dimensional simplices in $\Delta$ do not contain $\{u_1v_i,u_2v_i,\ldots,u_nv_i\}$ as a subset for any $i \geq 2$.  Conversely, if a $(k+1)$-element subset of $G$ does not contain $\{u_1v_i,u_2v_i,\ldots,u_nv_i\}$ as a subset for any $i \geq 2$, then there exists $m_i$ for each $i \geq 2$ such that $u_{m_i} v_i$ is not in the given subset, and this $(k+1)$-element subset is therefore a $k$-simplex that is a face of the maximal dimensional simplex $\mbox{conv} \{ G \backslash \{u_{m_1}v_2,u_{m_2}v_3,\ldots,u_{m_{h-1}}v_h\}\}.$  We conclude that a $(k+1)$-element subset of $G$ is a $k$-simplex in $\Delta$ if and only if it does not contain $\{u_1v_i,u_2v_i,\ldots,u_nv_i\}$ as a subset for any $i \geq 2$.  Thus, to determine $f_k$, we need to count the number of $(k+1)$-element subsets of $G$ that do not contain $\{u_1v_i,u_2v_i,\ldots,u_nv_i\}$ as a subset for any $i \geq 2$.  

Let us call a subset of the form $\{u_1v_i,u_2v_i,\ldots,u_nv_i\}$ a \emph{complete copy}. There are $\binom{(h+1)n}{k+1}$ $(k+1)$-element subsets of $G$, and the number of such subsets that contain $\ell$ complete copies as subsets is $\binom{hn - \ell n}{k+1 - \ell n} \binom{h-1}{\ell}$.  Thus by inclusion-exclusion,
\[
f_k = \sum_{\ell \geq 0} \binom{(h - \ell) n}{k+1 - \ell n} \binom{h-1}{\ell} {(-1)}^{\ell}.
\]
Since each maximal dimensional simplex in $\Delta$ has volume $\frac{1}{\dim(P(G))!}$, the result follows.
\qed

We now establish that Frobenius groups are two-level.  This relies on an important lemma in \cite{CollinsPerkinson2008}.

\begin{lemma}\label{frobeniusfacets}(See Corollary 4.5 in \cite{CollinsPerkinson2008})
The complement of any set of $|H|$ elements of $G$, one chosen from each of the cosets of $N$, forms the set of vertices of a facet of $P(G)$.  All facets of $P(G)$ arise this way.
\end{lemma}

\proof[Theorem~\ref{prp:frobeniustwolevel}] Let $J \subseteq G$ be the set of vertices of a facet of $G$.  Choose $H$ to be the set of coset representatives of $N$.  By Lemma~\ref{frobeniusfacets}, $J = G \backslash \{u_i, u_iv_1, u_iv_2, \ldots, u_iv_h\}$ for some fixed $i$. Now let \textbf{1} be the all ones vector in $\R^{h-1}$ and let $e_i$ be the standard basis vectors. Then we have
\[
\sum_{j \notin J} \overline{j} = \overline{u}_i + \overline{u_iv_1} + \overline{u_iv_2} + \cdots + \overline{u_iv_h} = \textbf{1} - e_1 - e_2 - \cdots - e_{k} = 0.
\]
Since $J$ was arbitrary, we conclude by Lemma~\ref{lem:galeexact} that $P(G)$ is two-level and thus $TH_1$-exact.
\qed

\section{Automorphism Groups of Binary Trees}~\label{AutBinTrees}

In this section, we present a method for computing the Ehrhart polynomials of groups that arise as automorphism groups of finite rooted binary trees.  The crux of this method lies in Theorem~\ref{thm:wreath}.  We first introduce some necessary group theoretic preliminaries.  For any groups $G \leq S_m$, $H \leq S_n$, the \emph{direct product} $G \times H \leq S_m \times S_m \leq S_{m+n}$ consists of elements $\{(g,h) : g \in G, \ h \in H \}$ with product $(g_1,h_1) \cdot (g_2,h_2) = (g_1g_2,h_1h_2)$.  By construction, the vertices of the permutation polytope of $G \times H$ are block matrices of the form $\{ g \oplus h : g \in G, \ h \in H\}$.  The \emph{wreath product} of $G$ by $S_n$, denoted $G \wr S_n$,  is the group $\{(g,h) : g \in G^{n}, h \in S_n\}$ under the operation defined by
$$(g',h') \cdot (g,h) = ((g'_1,g'_2,\ldots,g'_n),h') \cdot ((g_1,g_2,\ldots,g_n),h) := ((g'_{h'(1)}g_1,g'_{h'(2)}g_2,\ldots,g'_{h'(n)}g_n),h'h).$$
The vertices of the permutation polytope $P(G \wr S_n)$ are the $mn \times mn$ matrices $\{g \otimes h : g \in G, \ h \in S_n\}$.  For more on this, see \cite{Steinkamp1999}. 

We now prove that automorphism groups of rooted binary trees are always composed of direct products and wreath products of groups.  

\vspace{0.1in}

\begin{lemma}\label{lem:binaryautgroup}
Let $G$ be the automorphism group of a rooted binary tree.  Then $G$ can be written as a sequence of direct products of groups, and wreath products by symmetric groups of order at most two.
\end{lemma}

\proof
Label the vertices of $T$ by the positive integers $\{1,2,\ldots,n\}$ such that the root vertex is labeled $1$.  First assume the root of $T$ has one child, and without loss of generality assume its label is $2$.  Letting $T_2$ be the subtree of $T$ rooted at $2$, we have $\mbox{Aut}(T) = S_1 \times \mbox{Aut}(T_2)$.  Now assume instead that the root has two children that are labeled $2$ and $3$ without loss of generality.  Let $T_2$ be the subtree of $T$ rooted at $2$ and $T_3$ be the subtree of $T$ rooted at $3$.  If $T_2$ and $T_3$ are not isomorphic, then $\mbox{Aut}(T) = S_1 \times \left( \mbox{Aut}(T_2) \times \mbox{Aut}(T_3) \right)$.  If $T_2$ and $T_3$ are isomorphic, then $\mbox{Aut}(T) = S_1 \times \left( \mbox{Aut}(T_2) \wr S_2 \right)$.  The result then follows inductively.
\qed

The proof of Lemma~\ref{lem:binaryautgroup} indicates that computing the Ehrhart polynomial of groups arising as automorphism groups of rooted binary trees requires repeated computation of Ehrhart polynomials of direct products and wreath products by symmetric groups of order two.  Theorem~\ref{thm:wreath} indicates how Ehrhart polynomials behave under wreath products by symmetric groups of order two, and we prove this theorem now.

\vspace{0.1in}

\proof[Theorem~\ref{thm:wreath}] The vertices of the polytope $P(G \wr S_2)$ are precisely the matrices
\begin{eqnarray}\label{wreath}
\left \{
\begin{pmatrix}
X_1 & 0 \\
0 & X_2
\end{pmatrix},
\ \ 
\begin{pmatrix}
0 & X_1 \\
X_2 & 0
\end{pmatrix}\right\},
\end{eqnarray}
where the $X_i$ are vertices of $P(G)$.  Let $t \geq 2$ be an integer.  If $X_1,X_2$ are integer matrices in $kP(G)$ and $X_3,X_4$ are integer matrices in $(t-k)P(G)$, then
\begin{eqnarray}\label{wreathsum}
X=
\begin{pmatrix}
X_1 & 0 \\
0 & X_2
\end{pmatrix}
+
\begin{pmatrix}
0 & X_3 \\
X_4 & 0
\end{pmatrix}
\end{eqnarray}
is an integer matrix in $tP(G \wr S_2)$.  Moreover, if $X$ is an integer matrix in $tP(G \wr S_2)$, then there is a unique $k$ and unique integer matrices $X_1,X_2 \in kP(G)$ and $X_3,X_4 \in (t-k)P(G)$ such that $X$ can be decomposed as in (\ref{wreathsum}).  To see this, note that the supports of the above matrices imply that $X$ can be uniquely expressed as a convex combination of the $t^{th}$ dilations of matrices in $(\ref{wreath})$, so
\[
X=
\begin{pmatrix}
X_1 & 0 \\
0 & X_2
\end{pmatrix}
+
\begin{pmatrix}
0 & X_3 \\
X_4 & 0
\end{pmatrix}
\]
for some matrices $\{X_i \ | \ 1 \leq i \leq 4\}$.  Since the two summands have disjoint support, $X_1,X_2,X_3,X_4$ are all integer matrices.  Moreover, $X_1$ and $X_2$ have the same integer row and column sum, say $k$.  Consequently, $X_3$ and $X_4$ have row and column sum $(t-k)$.  Thus we conclude that the set of integer matrices in $tP(G \wr S_2)$ is in bijection with
\[
\bigcup_{k=0}^{n} \left( k(P(G) \times P(G)) \cap \Z^{n \times n} \right) \times \left( (t-k)(P(G) \times P(G)) \cap \Z^{n \times n} \right),
\]
and the result follows.
\qed


Theorem~\ref{thm:wreath} 
gives us a method for computing Ehrhart polynomials and hence volumes of permutation polytopes from groups arising as automorphism groups of rooted binary trees. First, given a rooted binary tree $T$, we compute the automorphism group Aut($T$) as a sequence of direct products and wreath products. Then we read the group Aut($T$) from left to right. If we encounter a direct product, we compute the Ehrhart polynomials of the corresponding groups and take the product of the polynomials. If we encounter a wreath product, we apply Theorem~\ref{thm:wreath}. This produces the Ehrhart polynomial of the permutation polytope associated to the tree $T$.


\begin{figure} [h]
\centering
\includegraphics[width=6cm]{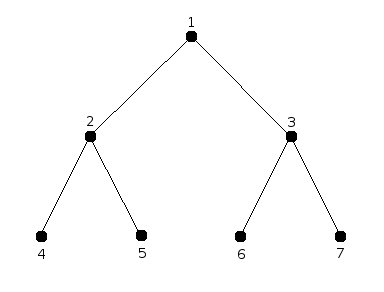}
\caption{A rooted binary tree T}
\label{fig:tree}
\end{figure}

\begin{example}\label{binautexample}

Consider the tree $T$ shown in Figure~\ref{fig:tree}.  Let $T_2,T_3$ be the subtrees rooted at $2$ and $3$ respectively.  Notice that
$\mbox{Aut}(T_2) = \mbox{Aut}(T_3)$ because $T$ is in fact unlabeled (we only place labels to illustrate how to compute the automorphism group).
The automorphism group of $T$ is therefore $S_1 \times \left[  \left( \mbox{Aut}(T_2) \right) \wr S_2 \right].$  Thus, by Theorem~\ref{thm:wreath}, its Ehrhart polynomial is 
\begin{align*}
i(P(Aut(T)),t) &= i(P(S_1),t) \cdot \left( \sum_{k=0}^t i^2(P(\mbox{Aut}(T_2)),k) \cdot i^2(P(\mbox{Aut}(T_2)),t-k) \right) \\
&= 1 \cdot \left(  \sum_{k=0}^t (k+1)^2 \cdot (t-k+1)^2 \right)  \\
&=  \sum_{k=0}^t (k+1)^2 \cdot (t-k+1)^2 \\
&= \sum_{k=0}^t k^4 + (-2t) \sum_{k=0}^t k^3 + (t^2 -2t -2) \sum_{k=0}^t k^2 + (t^2 -1) \sum_{k=0}^t k + (t+1)^2 \\
&=  \frac{1}{30}(t+1)(t+2)(t^2 + 4t + 5).
\end{align*}
By Theorem~\ref{thm:cyclicdihedralehrhart}, this is precisely the Ehrhart polynomial of $D_4$, which we should expect, since Aut$(T)$ is $S_1 \times D_4$ up to a relabeling of the generating set of $D_4$.  Moreover, we conclude that the volume of $P(Aut(T))$ is $\frac{1}{30}$.

\end{example}

We can further prove that for any rooted $T$, $P($Aut$(T))$ has a $P($Aut$(T))$-unimodular regular triangulation.

\begin{proposition}\label{prp:frobeniusunimodular}
If $T$ is a rooted tree, then $P(Aut(T))$ has a $P(Aut(T))$-unimodular regular triangulation.
\end{proposition}

\proof
Let $T$ be any rooted tree, and let $A_T$ be its adjacency matrix. Consider the polytope
\[
P_T = \left \{ X \in [0,1]^{n \times n} \ : \ A_T X = X A_T, \  \ \sum_{j=1}^n X_{ij} = 1 \ \forall i, \ \ \sum_{i=1}^n X_{ij} = 1 \ \forall j \right \}.
\]
The integer points of $P_T$ are permutations commuting with $A_T$, so they are precisely the automorphisms in Aut$(T)$.  Since $P_T$ is integral (see Theorem 2 of \cite{Tinhofer1986}), this implies that $P_T = P($Aut$(T))$.  But by Theorem 4.4 of \cite{DeLoeraHillarMalkinOmar2010}, $P_T$ is exact, which by Theorem 2.4 of \cite{Sullivant2006} and Theorem 4.2 of \cite{GouveiaParriloThomas2010} implies that any reverse lexicographic triangulation of $P($Aut$(T))$ is $P($Aut$(T))$-unimodular.  Since reverse lexicographic triangulations are regular, the result follows.
\qed

\section{Miscellaneous Permutation Polytopes}~\label{Miscellaneous}

In this section, we study miscellaneous permutation polytopes.  We begin by proving Proposition~\ref{prp:non-two-level}.  This proposition shows the difficulty of dealing with general permutation polytopes.

\vspace{0.1in}

\proof[Proposition~\ref{prp:non-two-level}] 
Since $P(A_2)$ and $P(A_3)$ have one and three vertices respectively, they are trivially two-level.  Since $A_4$ is a Frobenius group, Proposition~\ref{prp:frobeniustwolevel} implies that $P(A_4)$ is two-level.  For $n \geq 5$, by choosing $(\sigma,t,h) = (e,1,2)$ as in Theorem 3 of \cite{CunninghamWang2004}, we deduce that $P(A_n)$ has the facet defining inequality $\ell(x) \leq n-2$, where $\ell(x) = \sum_{j=3}^n x_{j,j} + \sum_{j=3}^n x_{j,1} + \sum_{j=3}^n x_{1,j}.$  Now $\ell(e) = n-2$, $\ell((1 \ 2)(4 \ 5))= n-4$, and $\ell((3 \ 4 \ 5)) = n-5$, and hence $P(A_n)$ is not two-level for $n \geq 5$.  To show that $P(A_n)$ is at least $(\lfloor \frac{n}{4} \rfloor +1)$-level for $n \geq 8$, we evaluate $\ell$ on $\sigma_i$, where $\sigma_0 = e$ and $\sigma_k = (1 \ 2)(3 \ 4) \cdots (4k-1 \ 4k)$ for $1 \leq k \leq \lfloor \frac{n}{4} \rfloor$.  
\qed

To conclude the paper, the table after the references lists the subgroups of $S_3$, $S_4$,
and $S_5$ and some their Ehrhart polynomials. Two groups stand out as incomplete, the alternating group $A_5$ and the
general affine group of degree one over the field of five
elements.  The latter group is generated by taking the semidirect product of
the additive and multiplicative groups of the field of five
elements, and is denoted by $GA(1,5)$.

\section*{Acknowledgements}

We would like to thank Igor Pak, David Perkinson, Raman Sanyal, Bernd
Sturmfels, and Rekha Thomas for helpful discussions. All authors were partially
supported by NSF grant DMS-0914107. The first and second authors were
supported by VIGRE NSF grant DMS-0636297. The third author was partially
supported by NSERC Postgraduate Scholarship 281174.

\bibliographystyle{plain}
\bibliography{biblio}

\vspace{0.1in}

\noindent\makebox[\textwidth]{
\begin{tabular}{ | c | c | c | p{3.5in} | }
\hline
\multicolumn{4}{ | c | }{Subgroups of $S_3$} \\ \hline
Order & Generators & Dim & Ehrhart Polynomial \\ \hline
$1$ & $<e> \cong \{e\}$ & 0 & $1$ \\ \hline
$2$ & $<(1 \ 2)> \cong C_2$ & 1 & $t + 1$ \\ \hline
$3$ & $<(1 \ 2 \ 3)> \cong C_3$ & 2 & $\frac{1}{2}t^2 + \frac{3}{2}t + 1$ \\ \hline
$6$ & $<(1 \ 2), (1 \ 3), (2 \ 3)> \cong S_3$ & 4 & $\frac{1}{8}t^4 + \frac{3}{4}t^3 + \frac{15}{8}t^2 + \frac{9}{4}t + 1$ \\ \hline


\multicolumn{4}{ | c | }{Subgroups of $S_4$} \\ \hline
Order & Group & Dim & Ehrhart Polynomial \\ \hline
$1$ & $<e> \cong \{e\}$ & $0$ & $1$ \\ \hline
$2$ & $<(1 \ 2)> \cong C_2$ & $1$ & $t + 1$ \\ \hline
$2$ & $<(1 \ 2)(3 \ 4)> \cong C_2$ & $1$ & $t + 1$ \\ \hline
$3$ & $<(1 \ 2 \ 3)> \cong C_3$ & $2$ & $\frac{1}{2}t^2 + \frac{3}{2}t + 1$ \\ \hline
$4$ & $<(1 \ 2), (3 \ 4)> \cong C_2 \times C_2$ & $2$ & $t^2 + 2t + 1$ \\ \hline
$4$ & \small{$<(1 \ 2)(3 \ 4), (1 \ 3)(2 \ 4)> \cong C_2 \times C_2$} & $3$ & $\frac{1}{6}t^3 + t^2 + \frac{11}{6}t + 1$ \\ \hline
$4$ & $<(1 \ 2 \ 3 \ 4)> \cong C_4$ & $3$ & $\frac{1}{6}t^3 + t^2 + \frac{11}{6}t + 1$ \\ \hline
$6$ & $<(1 \ 2), (1 \ 3), (2 \ 3)> \cong S_3$ & $4$ & $\frac{1}{8}t^4 + \frac{3}{4}t^3 + \frac{15}{8}t^2 + \frac{9}{4}t + 1$ \\ \hline
$8$ & $<(1 \ 2 \ 3 \ 4), (1 \ 2)(3 \ 4)> \cong D_4$ & $5$ & $\frac{1}{30}t^5 + \frac{1}{3}t^4 + \frac{4}{3}t^3 + \frac{8}{3}t^2 + \frac{79}{30}t + 1$ \\ \hline
$12$ & \tiny{$<(1 \ 2 \ 3), (1 \ 2 \ 4), (1 \ 3 \ 4), \ (2 \ 3 \ 4)> \cong A_4$} & $9$ & \tiny{$\frac{1}{5670}t^9 + \frac{1}{504}t^8 + \frac{23}{1890}t^7 + \frac{1}{15}t^6 + \frac{173}{540}t^5 + \frac{9}{8}t^4 + \frac{29797}{11340}t^3 + \frac{1199}{315}t^2 + \frac{383}{126}t + 1$} \\ \hline
$24$ & \tiny{$<(1 \ 2), (1 \ 3), (1 \ 4), (2 \ 3), (2 \ 4), (3 \ 4)> \cong S_4$} & $9$ & \tiny{$\frac{11}{11340}t^9 + \frac{11}{630}t^8 + \frac{19}{135}t^7 \frac{2}{3}t^6 + \frac{1109}{540}t^5 + \frac{43}{10}t^4 + \frac{35117}{5670}t^3 + \frac{379}{63}t^2 + \frac{65}{18}t + 1$} \\ \hline




\multicolumn{4}{ | c | }{Subgroups of $S_5$} \\ \hline
Order & Generators & Dim & Ehrhart Polynomial \\ \hline
$1$ & $<e> \cong \{e\}$ & 0 & $1$ \\ \hline
$2$ & $<(1 \ 2)> \cong C_2$ & 1 & $t + 1$ \\ \hline
$2$ & $<(1 \ 2)(3 \ 4)> \cong C_2$ & 1 & $t + 1$ \\ \hline
$3$ & $<(1 \ 2 \ 3)> \cong C_3$ & 2 & $\frac{1}{2}t^2 + \frac{3}{2}t + 1$ \\ \hline
$4$ & $<(1 \ 2), (3 \ 4)> \cong C_2 \times C_2$ & 2 & $t^2 + 2t + 1$ \\ \hline
$4$ & \small{$<(1 \ 2)(3 \ 4), (1 \ 3)(2 \ 4)> \cong C_2 \times C_2$} & 3 & $\frac{1}{6}t^3 + t^2 + \frac{11}{6}t + 1$ \\ \hline
$4$ & $<(1 \ 2 \ 3 \ 4)> \cong C_4$ & 3 & $\frac{1}{6}t^3 + t^2 + \frac{11}{6}t + 1$ \\ \hline
$5$ & $<(1 \ 2 \ 3 \ 4 \ 5)> \cong C_5$ & 4 & $\frac{1}{24}t^4 + \frac{5}{12}t^3 + \frac{35}{24}t^2 + \frac{25}{12}t + 1$ \\ \hline
$6$ & $<(1 \ 2 \ 3)(4 \ 5)> \cong C_6$ & 3 & $\frac{1}{2}t^3 + 2t^2 + \frac{5}{2}t + 1$ \\ \hline
$6$ & $<(1 \ 2), (2 \ 3), (1 \ 3)> \cong S_3$ & 4 & $\frac{1}{8}t^4 + \frac{3}{4}t^3 + \frac{15}{8}t^2 + \frac{9}{4}t + 1$ \\ \hline
$6$ & \footnotesize{$<(1 \ 2)(4 \ 5), (1 \ 3)(4 \ 5), (2 \ 3)(4 \ 5)> \cong S_3$} & 5 & $\frac{1}{40}t^5 + \frac{1}{8}t^4 + \frac{5}{8}t^3 + \frac{15}{8}t^2 + \frac{47}{20}t + 1$ \\ \hline
$8$ & \small{$<(1 \ 2 \ 3 \ 4), (1 \ 2)(3 \ 4)> \cong D_4$} & 5 & $\frac{1}{30}t^5 + \frac{1}{3}t^4 + \frac{4}{3}t^3 + \frac{8}{3}t^2 + \frac{79}{30}t + 1$ \\ \hline

$10$ & \small{$<(1 \ 2 \ 3 \ 4 \ 5), (2 \ 5)(3 \ 4)> \cong D_5$} & 8 & \tiny{$\frac{1}{8064}t^8 + \frac{5}{2016}t^7 + \frac{5}{192}t^6 + \frac{25}{144}t^5 + \frac{95}{128}t^4 + \frac{575}{288}t^3 + \frac{6515}{2016}t^2 + \frac{475}{168}t + 1$} \\ \hline
$12$ & \small{$<(1 \ 2 \ 3)(4 \ 5), (1 \ 2)(4 \ 5)> \cong D_6$} & 5 & $\frac{1}{8}t^5 + \frac{7}{8}t^4 + \frac{21}{8}t^3 + \frac{33}{8}t^2 + \frac{13}{4}t + 1$ \\ \hline
$12$ & \tiny{$<(1 \ 2 \ 3), (1 \ 2 \ 4), (1 \ 3 \ 4), \ (2 \ 3 \ 4)> \cong A_4$} & 9 & \tiny{$\frac{1}{5670}t^9 + \frac{1}{504}t^8 + \frac{23}{1890}t^7 + \frac{1}{15}t^6 + \frac{173}{540}t^5 + \frac{9}{8}t^4 + \frac{29797}{11340}t^3 + \frac{1199}{315}t^2 + \frac{383}{126}t + 1$} \\ \hline
$20$ & \footnotesize{$<(1 \ 2 \ 3 \ 4 \ 5), (1 \ 2 \ 4 \ 3)> \cong GA(1,5)$} & 9 & Too large to compute; volume=$\frac{19}{6538371840}$ \\  \hline
$24$ & \tiny{$<(1 \ 2), \ (1 \ 3), (1 \ 4), \ (2 \ 3), (2 \ 4), \ (3 \ 4)> \cong S_4$} & 9 & \tiny{$\frac{11}{11340}t^9 + \frac{11}{630}t^8 + \frac{19}{135}t^7 + \frac{2}{3}t^6 + \frac{1109}{540}t^5 + \frac{43}{10}t^4 + \frac{35117}{5670}t^3 + \frac{379}{63}t^2 + \frac{65}{18}t + 1$} \\ \hline
$60$ & \tiny{$<(1 \ 2 \ 3), (1 \ 2 \ 4), (1 \ 2 \ 5), (1 \ 3 \ 4), (1 \ 3 \ 5), (1 \ 4 \ 5),$} & 16 & Too large to compute \\
& \tiny{$(2 \ 3 \ 4), (2 \ 3 \ 5), (2 \ 4 \ 5), (3 \ 4 \ 5)> \cong A_5$} & & \\ \hline
$120$ & \tiny{$<(1 \ 2), (1 \ 3), (1 \ 4), (1 \ 5), (2 \ 3), (2 \ 4), (2 \ 5), (3 \ 4),$} & 16 & \tiny{$\frac{188723}{836911595520}t^{16} + \frac{188723}{20922789888}t^{15} + \frac{1008757}{5977939968}t^{14} + \frac{112655}{57480192}t^{13} + $} \\
& \tiny{$(3 \ 5), (4 \ 5)> \cong S_5$} & & \tiny{$\frac{72750523}{4598415360}t^{12} + \frac{984101}{10450944}t^{11} + \frac{125188639}{292626432}t^{10} + \frac{55426325}{36578304}t^9 + $}\\ 
& & & \tiny{$\frac{3541860299}{836075520}t^8 + \frac{196563587}{20901888}t^7 + \frac{3812839477}{229920768}t^6 + \frac{664118435}{28740096}t^5 + $} \\
& & & \tiny{$\frac{438177965089}{17435658240}t^4 + \frac{3028287247}{145297152}t^3 + \frac{6229735}{494208}t^2 + \frac{725}{144}t^1 + 1$} \\ \hline
\end{tabular} \label{tuttigrupi}
}

\end{document}